\documentclass[letterpaper, 10 pt, conference]{ieeeconf}  

\IEEEoverridecommandlockouts                              

\usepackage{tikz-cd}
\tikzcdset{/tikz/commutative diagrams/background color=red}
\tikzcdset{row sep/normal=1cm}

\usepackage{graphics} 
\usepackage{epsfig} 
\usepackage{amsmath} 
\usepackage{amssymb}  
\usepackage{booktabs}
\usepackage{makecell}
\usepackage{hyperref}

\newcommand{\dae}{f^{\mathrm{impl}}}
\newcommand{\ode}{f^{\mathrm{expl}}}
\newcommand{\integrator}{\Phi}

\newcommand{\tran}{^\top}

\newcommand{\R}{\mathbb{R}}
\newcommand{\nx}{{n_{\mathrm{x}}}}
\newcommand{\nuu}{{n_{\mathrm{u}}}}
\newcommand{\nz}{{n_{\mathrm{z}}}}
\newcommand{\dd}{{\mathrm{d}}}

\newcommand{\Tsim}{T_\mathrm{sim}}
\newcommand{\nstages}{s}
\newcommand{\nsteps}{n_{\mathrm{steps}}}

\title{\LARGE \bf
Fast integrators with sensitivity propagation for use in CasADi
}

\author{Jonathan Frey$^{1,2}$, Jochem De Schutter$^1$ and Moritz Diehl$^{1,2}$
\thanks{$^{1}$Department of Microsystems Engineering (IMTEK), University Freiburg, 79110 Freiburg, Germany
{\tt\small \{jonathan.frey, jochem.de.schutter, moritz.diehl\}@imtek.uni-freiburg.de}}%
\thanks{$^{2}$Department of Mathematics, University Freiburg, 79110 Freiburg, Germany
}%
}

\begin{document}

\maketitle
\thispagestyle{empty}
\pagestyle{empty}

\begin{abstract}
Efficient integrators with sensitivity propagation are an essential ingredient for the numerical solution of optimal control problems.
This paper gives an overview on the \texttt{acados} integrators, their \texttt{Python} interface and presents a workflow that allows using them with their sensitivities within a nonlinear programming (NLP) solver interfaced by \texttt{CasADi}.
The implementation is discussed, demonstrated and provided as open-source software.
The computation times of the proposed integrator and its sensitivity computation are compared to the native \texttt{CasADi} collocation integrator, \texttt{CVODES} and \texttt{IDAS} on different examples.
A speedup of one order of magnitude for simulation and of up to three orders of magnitude for the forward sensitivity propagation is shown for an airborne wind energy system model.
\end{abstract}

\section{Introduction}
Numerical optimization has become more widely used and tractable for increasingly complex problems due to advances in software and hardware.
When optimizing design and configurations of physical systems, one usually deals with optimization problems that are constrained by the evolution of a dynamic system model.
In the context of offline optimization, such models are often complex, stiff and associated with a high computational burden.
However, such accurate models can improve optimality compared to solutions obtained by less accurate models.

The software package \texttt{CasADi}~\cite{Andersson2019} is a widely used open-source framework for algorithmic differentiation and optimization which comes with interfaces to a variety of solvers for nonlinear programming (NLP) and initial value problems (IVP).
The open-source interior-point NLP solver \texttt{IPOPT}~\cite{Waechter2006} is probably the most widely used solver within \texttt{CasADi} and known for its numerical robustness.
However, an optimization solver derived from \texttt{CasADi} and \texttt{IPOPT} is typically not suitable for application in fast embedded optimal control applications.
In this context, more specialized software, like \texttt{acados}~\cite{Verschueren2021}, which implements SQP-type algorithms and integration methods that can efficiently propagate first and second order sensitivities, often needs to be used.
For efficiency, \texttt{acados} relies on the basic linear algebra package \texttt{BLASFEO}~\cite{Frison2020} and the (partial) condensing functionality implemented in \texttt{HPIPM}~\cite{Frison2020a}.

In this paper, we present a novel open-source software and workflow that allows one to use the \texttt{acados} integrators with their efficient native sensitivity propagation within a \texttt{CasADi} NLP obtained from a direct multiple shooting~\cite{Bock1984} parametrization.
The use of \texttt{acados} integrators within \texttt{IPOPT} allows solving NLPs which include the simulation of complex dynamical systems with one order of magnitude faster computation times compared to native collocation-based integrators provided by the \texttt{CasADi} software package.
Additionally, this workflow enables users to have a smoother transition from a \texttt{CasADi} + \texttt{IPOPT} implementation of their problem specific solver, to an \texttt{acados} implementation by using the intermediate step of using the \texttt{acados} integrator inside \texttt{IPOPT}.
\paragraph*{Outline}
The paper is organized as follows:
Section~\ref{sec:theory} gives a quick overview on integrators and sensitivity propagation for use within optimization.
Section~\ref{sec:implementation} gives an overview on the \texttt{acados} integrators, their interfaces to \texttt{Python}, and the \texttt{CasADi} wrapper.
Section~\ref{sec:experiments} shows some numerical experiments using the \texttt{acados} integrators within \texttt{CasADi}.

\section{Theoretical Background}
\label{sec:theory}
This section gives a brief overview on integrators for use within numerical optimization algorithms.\\
\subsection{Integrators for DAE and sensitivities}
Nonlinear dynamic systems are often given in form of a differential algebraic equation (DAE),
\begin{align}
\dae(\dot{x}(t), x(t), u(t), z(t)) = 0,
\label{eq:dae}
\end{align}
where the function $f^{\mathrm{impl}}: \R^\nx \times \R^\nx \times \R^\nuu \times \R^\nz \rightarrow \R^{\nx+\nz}$ describes the evolution of the dynamic system consisting of state $x$, control input $u$ and an algebraic state $z$ over time $t$.
We assume that the DAE \eqref{eq:dae} is of index-1, i.e. the matrix
$
\frac{\partial\dae}{\partial (\dot{x},z)}(\cdot)
$
is invertible and note that higher index DAEs are often reformulated as index-1 using index-reduction techniques.

We refer to an algorithm that solves the initial value problem (IVP) of DAE \eqref{eq:dae} together with an initial state $x_0 = x(0)$ and a given constant control input $u_0$ for a simulation time $\Tsim$ as an \textit{integrator}
\begin{align}
\integrator(x_0, u_0) \approx x(\Tsim).
\end{align}
Additionally, when solving optimization problems constrained by nonlinear dynamic systems, the integrator should be able to provide the forward sensitivities of the result with respect to the initial state and control input, i.e.
\begin{align}
\frac{\dd \integrator(x_0, u_0)}{\dd(x_0, u_0)} \in \R^{\nx\times(\nx+\nuu)}.
\end{align}
On the other hand, computing adjoint sensitivities directly is more efficient within some optimization algorithms, such as adjoint-based {inexact} SQP methods, or when BFGS hessian approximations are used.
The adjoint sensitivities for a given seed vector $\lambda \in \R^\nx$ are defined as
\begin{align}
\frac{\dd \lambda\tran\integrator(x_0, u_0)}{\dd(x_0, u_0)}\in \R^{1 \times (\nx+\nuu)}.
\end{align}
Moreover, second-order sensitivities are needed when using any optimization method that relies on an exact Hessian of the Lagrangian.
The second order sensitivities for a given seed vector $\lambda \in \R^\nx$ are given as
\begin{align}
\frac{\dd^2 \lambda\tran\integrator(x_0, u_0)}{\dd^2(x_0, u_0)} \in \R^{(\nx+\nuu)\times (\nx+\nuu)}.
\end{align}

\newcommand{\butcher}{{_{\mathrm{but}}}}
\subsection{Runge-Kutta methods in \texttt{acados}}
A Runge-Kutta method is defined by the equations
\begin{subequations}
\begin{align}
0 &= \dae(k_i, x_0 + \Tsim \sum_{j=1}^\nstages a_{ij}k_j, u_0, z_j) \label{eq:rk_eq} \\
 \nonumber &\quad\mathrm{for} \; i=1,\ldots,\nstages \\
x_+ &= x_0 + \Tsim \sum_{j=1}^\nstages b_j k_j. \label{eq:rk_out}
\end{align}
\end{subequations}
where the values $ a_{ij}, b_i, c_i $ for $ i,j=1,\ldots,\nstages$ form the Butcher tableau. 
In general, \eqref{eq:rk_eq} defines a set of implicit equations which typically are solved by applying a fixed number of Newton iterations.
Subsequently, the output is computed via \eqref{eq:rk_out}.

In case of an explicit ODE, of the form
\begin{align}
\dot{x}(t) = \ode(x(t), u(t)),
\end{align}
equation \eqref{eq:rk_eq} simplifies to a set of explicit equations, if the matrix $ A $ only has values below its diagonal.
Such methods are called explicit Runge-Kutta methods (ERK).

The more generally applicable implicit Runge-Kutta (IRK) methods have better stability properties and a higher order of integration.
The \texttt{acados} IRK method supports the Gauss-Legendre methods, which are of order $2s$ and A-stable, and the Gauss-Radau IIA methods, which are of order $2s-1$ and L-stable, which is desirable when handling stiff systems.
Additionally, it is common to use an equidistant partition of the simulation interval and apply the Runge-Kutta method $\nsteps$ times.

Within \texttt{acados}, ERK and IRK integrators with efficient sensitivity propagation (forward, adjoint and second-order) are implemented.
The \texttt{acados} IRK method implements the efficient symmetric Hessian propagation technique presented in~\cite{Quirynen2016}.
Interfaces to conveniently create integrators exist for \texttt{Python} and \textsc{Matlab}, which use the algorithmic differentiation (AD) and code generation functionality of \texttt{CasADi} to set up the ODE or DAE function and their derivatives.

When applying an SQP-type Optimal Control Problem (OCP) solver in an embedded context, it is typical to use an IRK method with a fixed number of Newton iterations~\cite{Quirynen2017} to reduce the variance of the OCP solution time towards a deterministic runtime.
However, outside of this context, it makes sense to (additionally) use a tolerance up to which the system of integration equations needs to be solved as a termination criterion, which has been added to the \texttt{acados} IRK implementation.

Additionally, the GNSF-IRK method~\cite{Frey2019} is implemented in \texttt{acados}, which can exploit different kinds of linear dependencies within a DAE and solve a nonlinear system in a reduced space instead of \eqref{eq:rk_eq}.
However, second-order sensitivity propagation is not implemented in this integrator yet, limiting its use to optimization methods without exact Hessians.

\subsection{Direct Multiple Shooting and Direct Collocation}
The two most general direct approaches for solving continuous time OCPs are (direct) multiple shooting and direct collocation.
Within multiple shooting, the dynamic system is discretized on each shooting interval using an integrator, which handles the integration variables ($k_i$ in the case of a RK method) internally.
On the other hand, a \textit{direct collocation} discretization directly handles the discretization within the optimization, by using the integration variables as optimization variables and the integration equations, \eqref{eq:rk_eq} in case of a RK method, as constraints of the NLP.
The direct collocation NLP has significantly more optimization variables and a more sparse structure.
\\
While there exist a variety of solvers that exploit the sparsity pattern of a multiple shooting based NLP or QP, solving the direct collocation based NLP is more attractive when using a solver that can exploit general sparsity patterns, such as \texttt{IPOPT}.
Moreover, an interior-point NLP solver can handle nonlinear dynamical systems via direct collocation more robustly, since the globalization strategies are applied on the whole system, including integration equations.
On the other hand, the multiple shooting based solution can often lead to a reduced memory footprint.
Additionally, it is essential for convergence on the multiple shooting based NLP that the sensitivity propagation of the integrator is consistent with its nominal evaluation, which can be achieved by freezing all adaptive settings (such as the internal step sizes, order, and iteration matrices).
%
%
The performance of the two discretization techniques highly depends on the dynamic system of interest, the required integration order and the available integrators and NLP solvers, which makes a general comparison very hard.

Note that in~\cite{Quirynen2017b} the concept of \textit{lifted integrators} was introduced which enables one to apply an SQP-type method with the classic OCP structured sparsity pattern, while iterating equivalently as on a direct collocation based problem.

\section{Implementation}
\label{sec:implementation}
This section gives an overview of existing integrators in \texttt{CasADi}, implementation details of the \texttt{acados} integrators that make them fast, and their \texttt{Python} interface.
Finally, the implementation of making the \texttt{acados} integrators fully usable in a \texttt{CasADi} NLP description and solver is discussed.
This feature has been released as open-source software~\cite{casados-integrators}.

\subsection{State-of-the-art in \texttt{CasADi}}
The \texttt{CasADi} \texttt{Integrator} class supports the following implementations:
\begin{itemize}
\item explicit RK-4 method
\item Collocation (IRK) method
\item \texttt{CVODES}~\cite{Serban2005}: ODE solver from the Sundials suite~\cite{Hindmarsh2005}
\item \texttt{IDAS}: DAE solver from the Sundials suite~\cite{Hindmarsh2005}
\end{itemize}
The \texttt{CasADi} collocation integrator can be used with Gauss-Legendre or Gauss-Radau IIA Butcher tableaux and internally uses the \texttt{Rootfinder} class.
In contrast to the IRK implementation in \texttt{acados}, it does not use the integration variables $k_i$ from equation \eqref{eq:rk_eq}, but instead the quantities $x_0 + \Tsim \sum_{j=1}^\nstages a_{ij}k_j$.

\texttt{CVODES} implements linear multistep methods such as Adams-Moulton (AM) methods and Backward differentiation formulas (BDF) with adaptive step-sizes.
\texttt{IDAS} mainly implements BDF formulas. 
Both \texttt{CVODES} and \texttt{IDAS} offer forward and adjoint sensitivity analysis.

The \texttt{CasADi} integrator interface allows one to specify quadrature states, which are quantities that do not influence the behavior of the dynamic system, but only depend on it and can thus be simulated more efficiently.
Note that this kind of quadrature states are a special case of a linear output system used in~\cite{Quirynen2013a} and the GNSF structure~\cite{Frey2019}.

Note that none of the baseline \texttt{CasADi} integrators can be code-generated, which could yield a speedup in run-time.


\subsection{\texttt{acados} integrators}
The model functions and their derivatives are code-generated by \texttt{CasADi}, then compiled and linked with the \texttt{acados} library, in the standard workflow of the high-level \texttt{acados} interfaces.
Thus, the \textit{creation time} of an \texttt{acados} integrator is significantly longer compared to the \texttt{CasADi} integrators mentioned above.
However, when using such an integrator inside an optimization solver, especially, when solving multiple problems, the reduced \textit{run-time} quickly makes up for the longer creation time.

All linear algebra operations in the implicit \texttt{acados} integrators are performed using \texttt{BLASFEO}.
The forward sensitivity propagation is performed efficiently using the implicit function theorem, following~\cite{Quirynen2017}.
The implicit \texttt{acados} integrators leave the integration variables from the last call in the memory, which has found to be useful in the context of MPC and leads to less Newton iterations in the numerical experiments in this paper.
It is possible to reset the initial guesses of the integration variables via the interfaces, which could be done in the proposed \texttt{CasADi} wrapper.

\subsection{\texttt{acados} integrator interface}

The \texttt{acados} interface to \texttt{Python} offers convenient methods to create an \texttt{acados} integrator, i.e., an instance of \texttt{AcadosSimSolver}, and to interact with it.
Possible interactions with \texttt{AcadosSimSolver} include setting the initial state, the control input, the parameter values, the guess (initialization) of the integration variables, the time horizon and the adjoint seed.
Most options, such as the size of the Butcher tableau $\nstages$, or the number of integration steps $\nsteps$ cannot be changed dynamically, since the memory size of the integrator depends on them and memory allocation should only happen once to create the integrator.
However, we added the possibility to set options for the sensitivity propagation modes, if the initial memory allocation took them into account.

The default \texttt{Python} wrapper of the \texttt{acados} solvers uses \texttt{ctypes}~\cite{ctypes}.
In order to speed up interactions with the \texttt{Python} object, we added a \texttt{cython} wrapper to the \texttt{AcadosSimSolver}.
This wrapper translates the code for the interaction layer into C code and compiles it~\cite{cython}.
Such a wrapper is also available for the \texttt{acados} OCP solver and in use in the open source driving assistance system \texttt{openpilot}~\cite{openpilot}.

\begin{figure}
    \begin{tikzcd}[description]
        \texttt{CasadosIntegrator}
        \arrow[d, "\texttt{get\_jacobian()}"]
        \arrow[dd, bend left, shift left=17ex, "\texttt{get\_reverse()}"] \\
        \texttt{CasadosIntegratorSensForw} \\
        \texttt{CasadosIntegratorSensAdj}
        \arrow[d, "\texttt{get\_jacobian()}"]
        \\
        \texttt{CasadosIntegratorSensHess} \\
    \end{tikzcd}
    \vspace{-.6cm}
    \caption{Derivative implementation using the \texttt{CasADi} \texttt{Callback}.}
    \vspace{-.4cm}
    \label{fig:casados}
\end{figure}
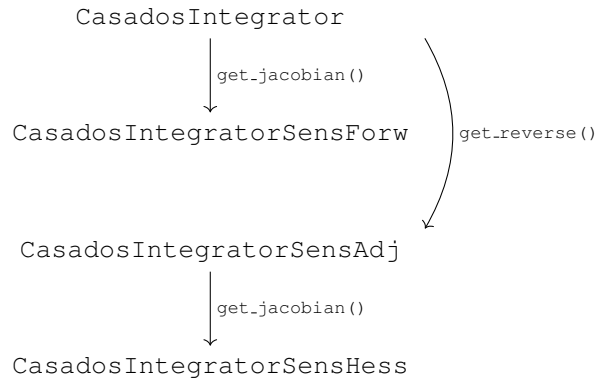

\subsection{\texttt{CasADi} wrapper for \texttt{acados} integrators}
We use the \texttt{CasADi} \texttt{Callback} class to wrap the interactions with the \texttt{acados} integrator into a \texttt{CasADi} \texttt{Python} object, called \texttt{CasadosIntegrator}.
The \texttt{CasADi} \texttt{Callback} class offers a public API to the \texttt{FunctionInternal} class.
We provide the implementation of the virtual methods that the \texttt{Callback} class needs to provide in order to be functional within the \texttt{CasADi} workflow.
Particularly, we implement forward and adjoint derivative computation via the functions \texttt{has\_jacobian()} and \texttt{get\_jacobian()}, respectively \texttt{has\_reverse()} and \texttt{get\_reverse()}, as illustrated in Figure~\ref{fig:casados}.
Internally, those derivative functions use the same \texttt{acados} integrator object.
This implementation is a major contribution of this paper and available in the open-source software package \texttt{casados-integrators} with the permissive 2-clause BSD license~\cite{casados-integrators}.

When the \texttt{CasadosIntegrator} or one of its derivatives is called, all the options of the \texttt{acados} integrator are set, such that the requested sensitivities, respectively the function values, are computed with minimal overhead.
When creating the \texttt{acados} integrator, all possible sensitivity options (forward, adjoint and Hessian) are activated, such that enough memory is allocated to execute all of them.
However, some of those options are temporarily disabled in the wrapper when they are not needed.

\subsection{Limitations of the \texttt{casados} integrators}
\label{sec:limitations}
When solving a multiple shooting OCP in \texttt{CasADi}, the integrator is called in the following order:
\begin{enumerate}
    \item nominal call for each shooting gap
    \item forward sensitivity call for each shooting gap
    \item adjoint call + Hessian call for each shooting gap
\end{enumerate}
This is suboptimal as the \texttt{acados} integrators offer to evaluate the nominal result and propagate sensitivity information within the same call, which allows to save some computations, such as the Newton iterations in case of an IRK.
Moreover, this separate evaluation contradicts the concept of freezing all adaptive components within the integrator.
Note that we observed no issues in this sense, in the numerical experiments, compared to other integrators, which have more adaptive components, such as \texttt{IDAS} and \texttt{CVODES}.

Additionally, the adjoint call before the Hessian call is unnecessary in case of the \texttt{acados}-based integrators.
This call stems from the fact that some functions need the non-differentiated output for the derivative evaluation\footnote{\url{https://github.com/casadi/casadi/issues/2019}}.

Possible extensions of this work include a similar Callback interface for \textsc{Matlab} and an implementation using Python memoryviews to reduce the overhead of the proposed Callback implementation.




\section{Numerical Experiments}
\label{sec:experiments}
The experiments presented here have been carried out using \texttt{acados} version v0.1.9 and \texttt{CasADi} 3.5.5 on a Lenovo T490s Laptop with Intel i5-8365U CPU, 16 GB RAM and Ubuntu 22.04.
We compare the computation times of the \texttt{casados} integrator with the ones provided by \texttt{CasADi} by solving an NLP with a nonlinear hanging chain model.
Second, we show how they can be used to solve a highly nonlinear NLP for an optimal orbit of an airborne wind energy system.
All timings are recorded using the \texttt{CasADi} internal timer, for the NLP solution in Section \ref{sec:chain} and the function evaluations in \ref{sec:awe} respectively.
The code for the benchmarks is publicly available at~\cite{casados-integrators}.

\subsection{Nonlinear chain of masses}
\label{sec:chain}
\begin{figure}[htb]
\includegraphics[width=\columnwidth]{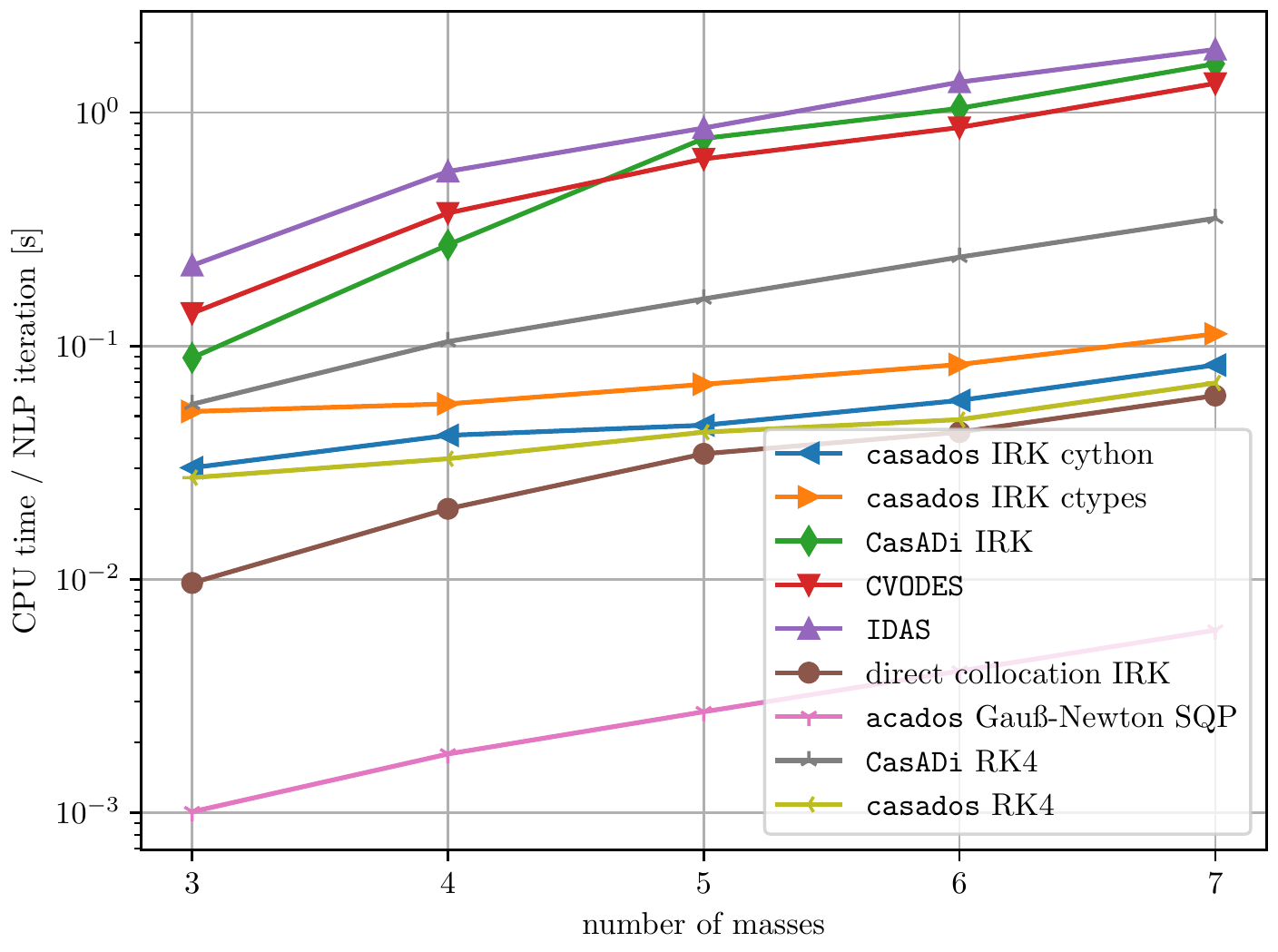}
\caption{CPU times with different integrators, IRK with $\nstages=2, \nsteps=1$.
}
\vspace{-.3cm}
\label{fig:chain_21}
\end{figure}
\begin{figure}[htb]
\includegraphics[width=\columnwidth]{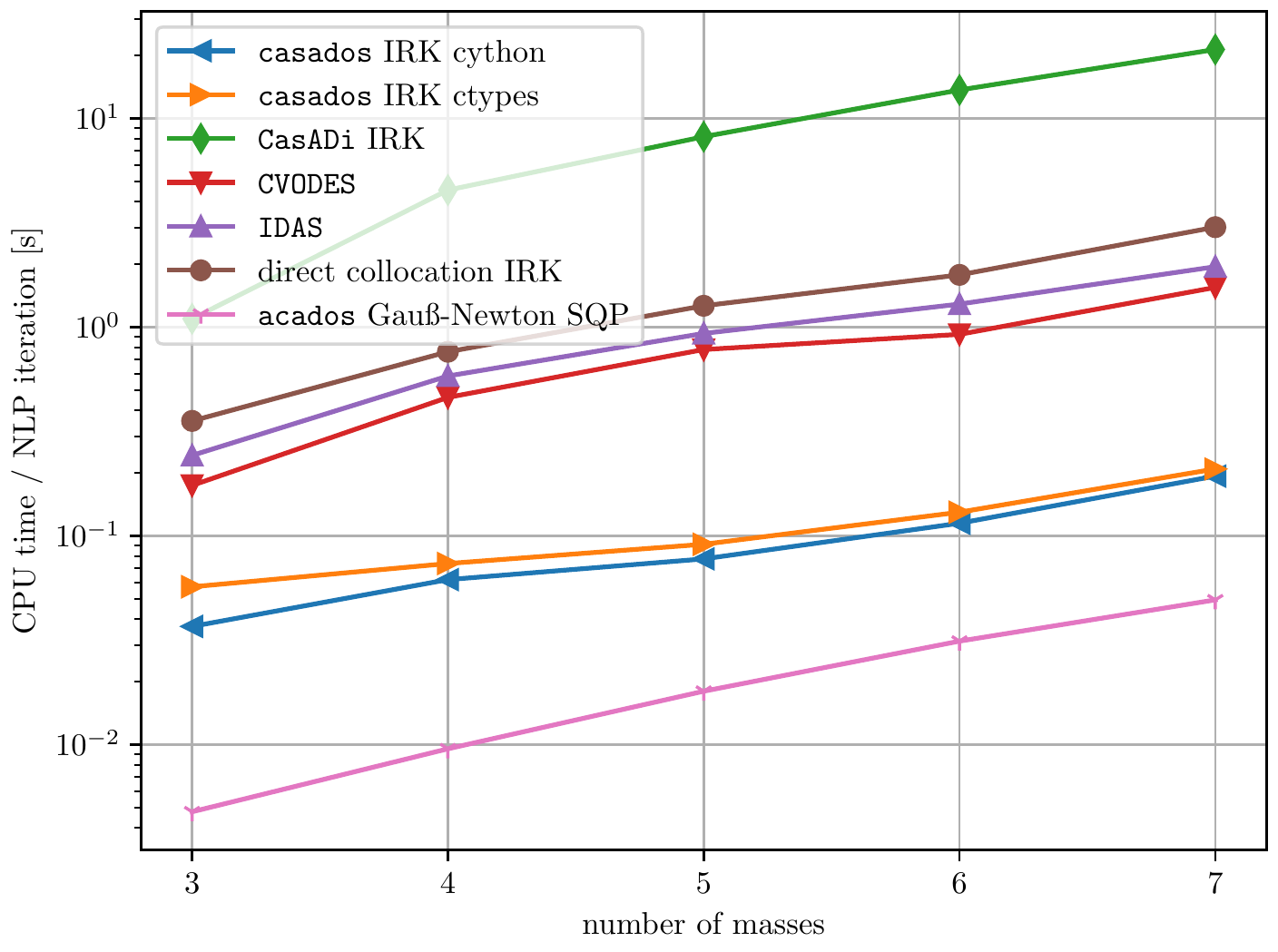}
\caption{CPU times with different integrators, IRK with $\nstages=4, \nsteps=4$.
}
\vspace{-.3cm}
\label{fig:chain_44}
\end{figure}

\newcommand{\nmass}{n_{\mathrm{mass}}}
We regard the problem of controlling a nonlinear hanging chain of masses used in~\cite{Zanelli2021} and originally proposed in~\cite{Wirsching2006} in order to see how the CPU time of an NLP solution scales with the state dimension for different integrators.
We can vary the number of masses in the chain and get a model with $\nx = 6(\nmass - 2)+3$.
We fix the horizon length to $N=40$.
The problem is solved using \texttt{IPOPT} and different integrators, namely \texttt{cacados} IRK and ERK with different wrappers and all available \texttt{CasADi} integrators, namely collocation (IRK), RK4 and CVODES, for the multiple shooting discretization.
Additionally, the same problem is solved with a direct collocation discretization and \texttt{IPOPT} and \texttt{acados} SQP with a Gau\ss-Newton Hessian approximation.
The settings of the different IRK version are equivalent and use the Gau\ss-Legendre methods.
We note that all integrator settings yield the same solution and the same number of NLP solver iterations.
The computation times are shown in Figure \ref{fig:chain_21} and~\ref{fig:chain_44}, for IRK  settings $\nstages=2, \nsteps =1$ and $\nstages=4, \nsteps =4$ respectively.
Although it is understood that the accuracy of RK4 is sufficient for the present example, we use the computationally expensive setting $\nstages=4, \nsteps =4$ to see how the computation times would evolve.
The timings are compared in terms of CPU time per NLP iterations in seconds, while the number of NLP solver iterations was between 11 and 15 for the \texttt{IPOPT} versions and between 5 and 6 for the \texttt{acados} SQP solution for varying $\nmass$.

Figure~\ref{fig:chain_21} gives a fair comparison of RK4 and the different IRK versions, since they are all of order $4$ and perform a single step.
We note, that even for small state dimensions, the corresponding \texttt{casados} integrator outperforms the native RK4 method and that the speedup grows significantly with the state dimension.
However, implicit methods still have better stability properties and are often preferred over explicit methods.
Note that in the \texttt{acados} implementation the gap between the ERK and IRK method is much smaller compared to the one between the corresponding \texttt{CasADi} integrators.
This can be attributed to the efficient \texttt{BLASFEO} linear algebra methods being used, as well as the tailored sensitivity propagation.

Comparing the differences between the two \texttt{casados} IRK integrators in Figure~\ref{fig:chain_21} and~\ref{fig:chain_44}, we see that for a computationally intensive integration scheme, the difference of the two versions is insignificant.
On the other hand, the speedup of the Cython wrapper is in the double digit percent range for simpler integration schemes, as in in Figure~\ref{fig:chain_21}.
Regarding Figure \ref{fig:chain_44}, we observe that the overall CPU time for the NLP solution with the equivalent IRK scheme is $\approx 100 $ times faster using the proposed integrator compared to the native \texttt{CasADi} IRK method.

Moreover, looking at the solution times of the equivalent direct collocation discretization, we observe that it is the fastest \texttt{IPOPT} variant overall for the simpler integration scheme in Figure \ref{fig:chain_21}, while being outperformed by a factor of 10 by the multiple shooting discretization using \texttt{casados} IRK in Figure \ref{fig:chain_44}.

The creation time, for a \texttt{casados} integrator with $\nmass=7$ is $3.1$s using \texttt{cython} and $0.8$s using \texttt{ctypes}, showing that its use can pay of after a single NLP iteration even for $\nstages=2, \nsteps =1 $.

The \texttt{acados} Gau\ss-Newton Hessian SQP algorithm consistently has the fastest iterations, due to multiple reasons.
The OCP-NLP structure is fully exploited, no second-order sensitivities are evaluated (Gau\ss-Newton), no Python overhead of the integrators, only a single call to the integrator for the nominal simulation and sensitivity propagation per shooting node gap and SQP iteration (compare to Sec.~\ref{sec:limitations}), fast QP solutions with \texttt{HPIPM}, \texttt{BLASFEO} and partial condensing~\cite{Frison2016} with a horizon of $10$.

\subsection{Airborne wind energy system}
\label{sec:awe}
\begin{figure}[tb]
\vspace{-.5cm}
	\input{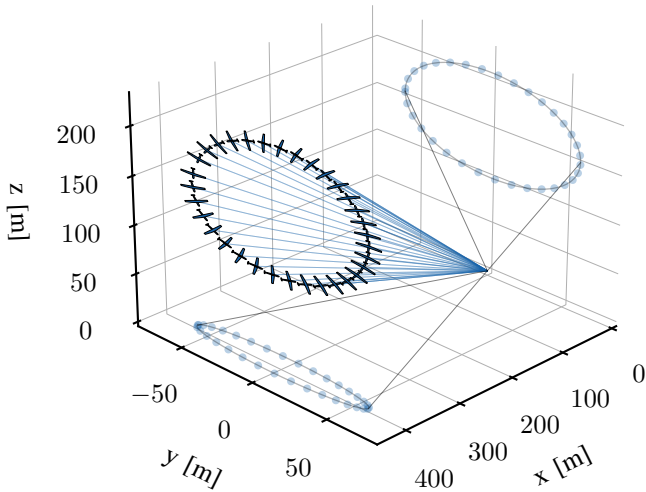}
	\caption{Optimal position and orientation trajectory of the drag-mode AWE system. The wind vector is pointing along the positive x-axis. To improve visibility, the aircraft size is enhanced with a factor 3 in the plot.}
	\label{fig:kite_traj}
\vspace{-.2cm}

\end{figure}

In the following, we consider the problem of finding a power-optimal periodic orbit for an airborne wind energy system.
This type of system consists of a tethered aircraft that flies fast crosswind loops in order to extract energy from the wind.
In this particular case, we consider a ``drag-mode" system~\cite{Vermillion2021}.
Here, electricity is generated on-board of the aircraft by means of small turbines that exert a braking force on the system, which is driven by the wind.
The electricity is transported to a ground station through the tether.

For this study, we choose the reference model presented in~\cite{Malz2019}.
For the on-board turbine drag force and power, we use the expressions given in~\cite{Zanon2013a}.
The model describes the full six-degree-of-freedom dynamics of a small AWE system with 5.5 m wing span aircraft.
The system state is modeled in non-minimal coordinates so as to tailor the dynamic equations for use in Newton-type optimization algorithms.
Hence, the resulting DAE dynamics contain six invariants, which are stabilized using a Baumgarte scheme~\cite{Gros2018}.
The system dimensions are $n_{\mathrm{x}}~=~23, n_{\mathrm{u}}~=~4, n_{\mathrm{z}}~=~1$.

The dynamics are discretized using an IRK method with Gau\ss-Radau IIA Butcher tableau and $\nstages=4$, with sampling time  $\Tsim = 0.3364 \ \mathrm{s}$.
These discrete dynamics are used to formulate the periodic discrete-time optimal control problem as stated in~\cite[Eq. 3]{DeSchutter2020}, with $N = 40$ intervals.
We impose the realistic flight envelope constraints given in~\cite{Licitra2019a}, resulting in 21 linear and 9 nonlinear inequality constraints per interval.
The discrete stage cost is defined as the negative average on-board energy generated in one interval, combined with the regularization terms given in~\cite{Licitra2019a}.

The resulting NLP is highly non-convex and requires a good initial guess to converge.
Therefore we solve the problem based on an initial guess provided by the open-source AWE optimization toolbox~\texttt{awebox}~\cite{awebox2020}.
The power-optimal flight trajectory is visualized in Figure~\ref{fig:kite_traj}.

\paragraph{NLP solution timing comparison}
Table \ref{tab:time_kite_dc_vs_ms} shows the NLP solution times for direct collocation and multiple shooting with \texttt{casados} IRK respectively, computed using \texttt{IPOPT} and \texttt{MA57}.
Both methods result in timings in the same order of magnitude, with direct collocation outperforming multiple shooting by a factor of 2.
While the multiple shooting timings are dominated by the NLP function evaluations, the step computation timings are reduced by one order of magnitude compared to direct collocation.
This implies that the potential speedup of a parallelized linearization is significantly higher for the multiple shooting discretization.

\begin{table}[]
	\centering
	\begin{tabular}{ccc}
		& \thead{multiple shooting \\(\texttt{casados} IRK)} & direct collocation \\ \toprule
		total & 4.12 & 2.02  \\
		NLP function eval. & 3.94 & 0.537 \\
		step computation & 0.174 & 1.485 \\
                NLP solver iterations & 20 & 10
	\end{tabular}
	\caption{AWE system example: NLP solution timings in s.}
        \vspace{-.5cm}
	\label{tab:time_kite_dc_vs_ms}
\end{table}

\paragraph{Integrator Timing Comparison}
Since we did not manage to solve the NLP described in the previous paragraph with the other integrators mentioned here despite trying various options, we limit the comparison of computation times to forward simulation and sensitivity propagation in this paragraph.
In Table~\ref{tab:time_kite_sim}, we compare the computation times of forward simulating the states and controls from the optimal trajectory in Figure~\ref{fig:kite_traj}.
We observe that a speedup factor larger than 30 is achieved consistently.

\begin{table}[]
\centering
\begin{tabular}{cccc}
        & \texttt{casados} & \texttt{CasADi} IRK & speedup factor \\ \toprule
        median & 0.2554 & 7.0845 & 27.74 \\
        max & 0.4911 & 8.4917 & 17.29 \\
        min & 0.2283 & 6.4680 & 28.33
\end{tabular}
\caption{Timings in ms, forward simulation of given initial state and control trajectory.}
\vspace{-.5cm}
\label{tab:time_kite_sim}
\end{table}

In Table~\ref{tab:time_kite_jac}, we compare the computation times of computing the Jacobian of the solution w.r.t. initial state and control input along the same trajectory.
A comparison of those Jacobians shows that the maximum difference between any corresponding entries is below $10^{-12}$, thereby confirming that the methods are mathematically equivalent.
We note that for this kind of system the proposed integrator gives a speedup of roughly a factor $ 1000$.

\begin{table}[]
\centering
\begin{tabular}{cccc}
        & \texttt{casados} & \texttt{CasADi} IRK & speedup factor \\ \toprule
        median & 0.6839 & 789.6 & 1154.6 \\
        max & 0.9917 & 884.4 & 891.8 \\
        min & 0.6528 & 744.4 & 1140.4
\end{tabular}
\caption{Timings to evaluate Jacobian of simulation result w.r.t. initial state and control input for different initial states and control inputs.}
\vspace{-.5cm}
\label{tab:time_kite_jac}
\end{table}




\section{Conclusions \& Outlook}
In this paper, we gave a detailed overview on the \texttt{acados} integrators and their interfaces to Python.
We introduced a novel interface that makes those integrators available in \texttt{CasADi} Python via the open-source \texttt{casados} integrator software package.
We show speedups of a factor of up to 1000 with respect to the integrators distributed by CasADi, when using them in an NLP solver.
This enables the solution of problems that cannot be solved with the standard integrators in reasonable computation times.







\section*{Acknowledgments}
This research was supported by DFG via Research Unit FOR 2401 and project 424107692 and by the EU via ELO-X 953348.
We thank Katrin Baumgärtner for fruitful discussions, Mario Zanon for asking for such a feature in the \texttt{acados} forum and Joris Gillis for advice regarding the Callback implementation.

\bibliographystyle{ieeetr}
\bibliography{syscop}
\addcontentsline{toc}{chapter}{Bibliography}

\end{document}